\documentclass[12pt,a4paper]{article}
\usepackage{graphicx}
\usepackage{amsfonts}
\usepackage{amsmath}
\usepackage{amssymb}
\usepackage[ansinew]{inputenc}
\usepackage{latexsym}
\usepackage{tabularx}
\usepackage{makeidx}
\usepackage{color}
\allowdisplaybreaks
\newtheorem{thm}{Theorem}

\newtheorem{lem}{Lemma}

\newenvironment{proof}{\begin{trivlist}
\item[\hskip\labelsep{\it Proof.}]}{$\hfill\Box$\end{trivlist}}

\newcommand{\eps}{\varepsilon}

\def\({\left(}
\def\){\right)}
\def\fl#1{\left\lfloor#1\right\rfloor}

\makeatletter
\newcommand{\rdots}{\mathinner{\mkern1mu\lower-1\p@\vbox{\kern7\p@\hbox{.}}
\mkern2mu \raise4\p@\hbox{.}\mkern2mu\raise7\p@\hbox{.}\mkern1mu}} \makeatother

\begin{document}

\title{On the Distribution of the Euler Function of Shifted Smooth Numbers}

\author{
Stefanie S. Loiperdinger\thanks{Supported by an International Macquarie University Research Excellence Scholarship} \ and
Igor E. Shparlinski\thanks{Supported by   ARC grant DP0556431}\\
\\
Department of Computing\\
Macquarie University\\
Sydney, NSW, 2109\\
Australia\\
\tt{\{stefanie,igor\}@ics.mq.edu.au} }

\date{\today}
\maketitle

\begin{abstract}
We give asymptotic formulas for some average values of the
Euler function on shifted smooth numbers. The result is based
on various estimates on the distribution of
smooth numbers in arithmetic progressions which are due
to A.~Granville and {\'E}.~Fouvry \& G.~Tenenbaum.
\end{abstract}

\section{Introduction}

 An integer $n \geq 1$ is called $y$-smooth if every prime factor $p$ of $n$
 satisfies $p\leq y$. For a detailed introduction to smooth numbers, their properties and applications, see~\cite{CEP,Granville93a, Granville93b, GranvilleXX, HildTenen93,Sound,Ten} and references therein.

 We denote by $S(x,y)$ the set of numbers less than or equal to
$x$ that are $y$-smooth, that is,
$$S(x,y)=  \{n:1 \leq n \le x \text{ and } n \text{ is }
y\text{-smooth}  \}. $$ Furthermore let $\Psi (x,y) = \# S(x,y)$ be the counting function for smooth numbers.

Also, as usual, we use $\varphi(k)$ to denote the Euler function of an integer $k \ge 1$.

In this paper, we obtain asymptotic formulas for some  average values of the Euler function of shifted smooth numbers. Namely, for real $x \ge
y \ge 2$, we define
$$T(x,y)=\sum_{\substack{a<n \le x\\
n\in S(x,y)} }\frac{\varphi(n - a)}{n - a} \qquad\text{and}\qquad
V(x,y)=\frac{1}{\Psi(x,y)}\sum_{\substack{a<n \le x\\
n\in S(x,y)} } \varphi(n - a) ,
$$
where $a\ne 0 $ is a fixed  integer
(throughout the paper, the implied constant may depend on $a$).

\section{Preparations}

Throughout the paper, we use
$U = O(V)$, $U \ll V$, and $V \gg U$
as  equivalents of the inequality $|U| \le c V$ with some constant
$c> 0$, which  may depend only on  $n$.

We recall that the {\it Dickman--de~Bruijn\/}
function $\rho(u)$ is defined by
$$
\rho(u) = 1, \qquad 0 \le u \le 1,
$$
and
$$
\rho(u) = 1 - \int_{1}^u \frac{\rho(v- 1)}{v}\,d v, \qquad u  > 1.
$$
Then, by~\cite[Chapter~III.5, Corollary~9.3]{Ten}, we have

\begin{lem}\label{psi_rho}
For any $\eps>0$, the estimate
$$\Psi(x,y)=  x \rho(u) \(1 + O\(\frac{ \log (u+1)}{\log y} \)\) $$
holds uniformly in the range
$$\exp((\log \log x)^{5/3+\eps})\leq y \leq x,
$$
where
$$
u=\frac{\log x}{\log y}.
$$
\end{lem}

The following asymptotic estimate on $\rho(u)$ follows immediately from
a much more precise result of~\cite[Chapter~III.5, Theorem~8]{Ten}.

\begin{lem}\label{rho}
For any $u \to \infty $, we have
$$\rho(u) = \exp\(-(1 + o(1))u \log u\).
$$
\end{lem}

We note that the bound
\begin{equation}
\label{psi_CEP}
\Psi(x,y) = xu^{-u+o(u)},
\end{equation}
due to Canfield, Erd{\H o}s and Pomerance~\cite[Corollary to
Theorem~3.1]{CEP},
holds in a much wider range than one can obtain from
Lemmas~\ref{psi_rho} and~\ref{rho},
see also~\cite{GranvilleXX,HildTenen93,Ten}.

Furthermore, the following upper bound on the
derivative of $\rho(u)$ is a very weak form of a much more precise result~\cite[Chapter~III.5,
Corollay~8.3]{Ten}.

\begin{lem}\label{rho'}
For any $u >0$, we have
$$\rho'(u) \ll \rho(u) \log(u+1) .
$$
\end{lem}

For any integers $a$ and $d$ with gcd$(a,d)=1$, let
$$\Psi(x,y;a,d)=\#\{n\in S(x,y):n\equiv a\pmod d\} $$ and let
$$\Psi_d(x,y)=\#\{n\in S(x,y): \gcd(n,d)=1\}. $$
In general, one expects that
$$
\Psi(x,y;a,d)\sim \frac{\Psi_d(x,y)}{\varphi(d)},
$$
for sufficiently large $x$.
%

Granville~\cite{Granville93a} has proved the following bounds on the average of smooth numbers lying in a fixed arithmetic progression.

\begin{lem}\label{psi_d_a}
Let $A$ be a fixed positive number. Then there exist positive constants $\gamma$ and $\delta$, depending only on $A$,
 such that for
$$
\Delta = \min\left\{\exp\(\gamma \frac{\log y \log \log y}{\log \log \log y}\),\frac{\sqrt{x}}{(\log x)^\delta}\right\}
$$
uniformly over   $y \geq 100$ we have
$$\sum_{d\leq \Delta} \max_{z\leq x}\max_{\gcd(a,d)=1}\left| \Psi(z,y;a,d)-\frac{\Psi_d(z,y)}{\varphi (d)} \right|=O\(\frac{\Psi(x,y)}{(\log y)^A}\), $$
where the implied constant depends only on $A$.
\end{lem}

Finally, Fouvry and Tenenbaum~\cite{FouvryTenen} give the
following asymptotic formula for the number of smooth numbers that are coprime to
$d$.

\begin{lem}\label{psi_d}
For any $\eps>0$ there exist $x_0(\eps)$ such that for $x\ge x_0(\eps)$, the estimate
$$\Psi_d(x,y)= \frac{\varphi(d)}{d}\Psi(x,y)\(1+O\(\frac{\log \log (dy) \log \log x}{\log y} \)\) $$
holds uniformly in the range
$$\exp((\log \log x)^{5/3+\eps})\leq y \leq x, \qquad \log \log(d+2)\leq \(\frac{\log
y}{\log(u+1)}\) ^{1-\eps},$$ where
$$
u=\frac{\log x}{\log y}.
$$
\end{lem}

\section{Asymptotic Formulas}

We are now ready to obtain our main results.

\begin{thm}\label{thm:Txy}
There exists an absolute constant $C > 0$ such that for a sufficiently large $x$ the bound
$$
T(x,y)= \Psi(x,y)
 \(\frac{6}{\pi^2} + O\(\frac{\log \log x \log \log y}{\log y}\)\)
$$
holds uniformly  in the range
$$
x \ge y \ge  \exp\(C \sqrt{\log x \log \log \log x}\).
$$
\end{thm}

\begin{proof}
Using the well known identity
$$\varphi(n)=n\sum_{d|n}\frac{\mu(d)}{d},$$
where $\mu(d)$ is the M\"{o}bius function, see~\cite[Equation~(16.3.1)]{HardyWright}, and changing the order of summation, we can rewrite
$T(x,y)$ in the following way,
\begin{eqnarray*}
\lefteqn{T(x,y)
= \sum_{\substack{a<n \le x\\
n\in S(x,y)} } \sum_{d|(n - a)}\frac{\mu(d)}{d}}\\
&& \qquad\qquad
 = \sum_{d\leq x}\frac{\mu(d)}{d}\sum_{\substack{a<n \le x \\
n\in S(x,y)\\ n\equiv a\pmod d}}1= \sum_{d\leq x}\frac{\mu(d)}{d}\Psi(x,y;a,d).\\
\end{eqnarray*}
Let $\gamma$ and $\delta$ are chosen as in Lemma~\ref{psi_d_a},
corresponding to $A=1$.
We now  define
$$\Delta=\min\left\{\exp\(\gamma \frac{\log y \log \log y}{\log \log \log y}\),
\frac{\sqrt{x/a}}{(\log x/a)^\delta}\right\},
$$
and write
\begin{equation}
\label{eq:split} T(x,y) = \sum_{d\leq x}\frac{\mu(d)}{d}\Psi(x,y;a,d) = \Sigma_1 + \Sigma_2 ,
\end{equation}
where
\begin{eqnarray*}
\Sigma_1 &=& \sum_{d\leq \Delta}\frac{\mu(d)}{d}\Psi(x,y;a,d);\\
\Sigma_2 &=& \sum_{x \ge d > \Delta}\frac{\mu(d)}{d}\Psi(x,y;a,d).
\end{eqnarray*}

For  $\Sigma_1$ we have
\begin{equation}
\label{eq:Sigma1}
\Sigma_1 =\sum_{d \le  \Delta} \frac{\mu(d) \Psi_d(x,y)}{d\varphi(d)} + O(R),
\end{equation}
where
$$
R = \sum_{d \le \Delta} \frac{1}{d}
\left|\Psi(x,y;a,d) -\frac{\Psi_d(x,y)}{\varphi(d)} \right|.
$$
Now, for each divisor $f\mid a$, we collect together the terms with
$\gcd(a,d) = f$, getting
\begin{equation}
\label{eq:Rf}
R=\sum_{f \mid a} R_f,
\end{equation}
where
\begin{equation}
\label{eq:Rf Psi-df}
\begin{split}
R_f &= \sum_{\substack{d \le \Delta \\
\gcd(a,d)=f}} \frac{1}{d} \left|\Psi(x,y;a,d) -\frac{\Psi_d(x,y)}{\varphi(d)} \right|\\
 &= \sum_{\substack{d \le \Delta \\
\gcd(a,d)=f}} \frac{1}{d} \left|\Psi(x/f,y;a/f,d/f) -\frac{\Psi_d(x,y)}{\varphi(d)} \right|,
\end{split}
\end{equation}
provided that $y > |a|$.
We now note that
Lemma~\ref{psi_d} implies
that
\begin{equation}
\label{eq:Psi-df}
\begin{split}
\frac{\Psi_{d/f}(x/f,y)}{\varphi(d/f)} & =
 \frac{1}{d/f}\Psi(x/f,y)\(1+O\(\frac{\log \log (dy) \log \log x}{\log y} \)\) \\
 & =
 \frac{f}{d}\Psi(x/f,y)\(1+O\(\frac{\log \log (dy) \log \log x}{\log y} \)\) .
\end{split}
\end{equation}
Furthermore, denoting
$$
u_f = \frac{\log (x/f)}{\log y} = u + O\(\frac{1}{\log y} \),
$$
we see from Lemma~\ref{rho'} that
$$
\rho(u_f) = \rho(u) \(1+O\(\frac{\log (u+1)}{\log y} \)\).
$$
Thus, by  Lemma~\ref{psi_rho} we have
$$
\Psi(x/f,y) = \frac{1}{f} \Psi(x,y)\(1+O\(\frac{\log (u+1)}{\log y} \)\).
$$
 Therefore~\eqref{eq:Psi-df} can be re-written as
\begin{eqnarray*}
\frac{\Psi_{d/f}(x/f,y)}{\varphi(d/f)} &= &
 \frac{1}{d}\Psi(x,y)\(1+O\(\frac{\log \log (dy) \log \log x}{\log y}+\frac{\log (u+1)}{\log y} \)\) \\
&= &
 \frac{1}{d}\Psi(x,y)\(1+O\(\frac{\log \log (dy) \log \log x}{\log y} \)\)
\end{eqnarray*}
(since $u \ll \log x$).
Applying Lemma~\ref{psi_d} again, we obtain
$$
\frac{\Psi_d(x,y)}{\varphi(d)}
= \frac{\Psi_{d/f}(x/f,y)}{\varphi(d/f)} \(1+O\(\frac{\log \log (dy) \log \log x}{\log y} \)\).
$$
 Accordingly, since the series
$$
 \sum_{\substack{d =1 \\
\gcd(a,d)=f}}^\infty \frac{1}{d\varphi(d/f)}  < \infty
$$
converges,
we now derive from~\eqref{eq:Rf Psi-df} that
\begin{eqnarray*}
R_f & = &\sum_{\substack{d \le \Delta \\
\gcd(a,d)=f}} \frac{1}{d} \left|\Psi(x/f,y;a/f,d/f) -
\frac{\Psi_{d/f}(x/f,y)}{\varphi(d/f)}  \right|\\
& & \qquad \qquad \qquad\qquad
+~O\(\Psi(x,y) \frac{ \log \log (\Delta y) \log \log x}{\log y} \)\\
& \ll &\sum_{\substack{d \le \Delta \\
\gcd(a,d)=f}}  \left|\Psi(x/f,y;a/f,d/f) -
\frac{\Psi_{d/f}(x/f,y)}{\varphi(d/f)}  \right|\\
& & \qquad \qquad \qquad \qquad
+~\Psi(x,y) \frac{ \log \log (\Delta y) \log \log x}{\log y} .
\end{eqnarray*}
Moreover, in the considered range of $x$ and $y$,
for sufficiently large $x$, we have
$$
y \le \Delta^3,
$$
hence
$$
\log \log (\Delta y) \le \log \log (\Delta^4) \le \log \(\frac{4\gamma \log y \log \log y}{\log \log \log y}\) = O(\log \log y).
$$

Since $\gamma$ and $\delta$ in the definition of $\Delta$
are chosen to correspond to  $A=1$  in Lemma~\ref{psi_d_a},
we obtain
$$
R_f  \ll \Psi(x,y) \frac{ \log \log x \log \log y}{\log y}
$$
which after the substitution into~\eqref{eq:Rf} yields
\begin{equation}
\label{eq:R bound}
R  \ll \Psi(x,y) \frac{ \log \log x \log \log y}{\log y} .
\end{equation}

 We see that for
$d \le \Delta$ the condition of Lemma~\ref{psi_d}
$$\log \log (d+2) \leq \(\frac{\log y}{\log (u+1)}\)^{1-\eps}
$$
is satisfied (provided $x$ is large enough), so we derive
\begin{equation}
\label{eq:asymp1}
\begin{split}
\sum_{d \le  \Delta} \frac{\mu(d) \Psi_d(x,y)}{d\varphi(d)}
& =~ \Psi(x,y)\sum_{d\le  \Delta}\frac{\mu(d)}{d^2}\(1+O\(\frac{\log \log (dy) \log \log x}{\log y} \)\)\\
 =~&\Psi   (x,y) \(\sum_{d\le  \Delta}\frac{\mu(d)}{d^2}  +O\(
 \frac{ \log \log x}{\log y}
\sum_{d\le  \Delta}  \frac{\log \log (dy)  }{d^2} \)\).
\end{split}
\end{equation}
We also have
\begin{equation}
\label{eq:main term} \sum_{d \leq  \Delta} \frac{\mu(d)}{d^{2}}= \sum_{d=1}^\infty \frac{\mu(d)}{d^{2}}+ O\(\sum_{d \geq  \Delta}
\frac{1}{d^{2}}\)=\frac{1}{\zeta (2)}+O\(\frac{1}{ \Delta}\)= \frac{6}{\pi^2}+O\(\frac{1}{ \Delta}\),
\end{equation}
where $\zeta(s)$ is the Riemann zeta-function, see~\cite[Theorem~287 and Equation~(17.2.2)]{HardyWright}. To estimate the error term
in~\eqref{eq:asymp1} we use the trivial inequality
\begin{equation}
\label{eq:error term} \sum_{d\le  \Delta}  \frac{\log \log (dy)  }{d^2} \le \sum_{d\le  \Delta}  \frac{\log \log ( \Delta y) }{d^2}   \ll \log
\log ( \Delta y) \ll \log \log y.
\end{equation}
Thus, substituting~\eqref{eq:main term} and~\eqref{eq:error term} in~\eqref{eq:asymp1}, we derive
\begin{equation}
\label{eq:prelim} \sum_{d \le  \Delta}  \frac{\mu(d) \Psi_d(x,y)}{d\varphi(d)}=\Psi(x,y)\(\frac{6}{\pi^2}+O\(\frac{\log \log x \log \log
y}{\log y}\)\).
\end{equation}

Combining~\eqref{eq:R bound} and~\eqref{eq:prelim}, we
deduce from~\eqref{eq:Sigma1}
\begin{equation}
\label{eq:Sigma1 asymp} \Sigma_1  =
\Psi(x,y)\(\frac{6}{\pi^2}+O\(\frac{\log \log x \log \log y}{\log y}\)\).
\end{equation}

For $\Sigma_2$ we have the trivial estimate
\begin{equation}
\label{eq:Sigma2 bound}
\begin{split}
|\Sigma_2| & \le~ \sum_{x \ge d > \Delta}\frac{1}{d} \sum_{\substack{a<n \le x \\ n\in S(x,y)\\ n\equiv a\pmod d}}1  \le   \sum_{x \ge d >
\Delta}\frac{1}{d}
\sum_{\substack{a<n \le x \\ n\equiv a\pmod d}} 1 \\
 & \le~\sum_{x \ge d > \Delta}\frac{1}{d} \(\fl{x/d} + 1\)
 \le 2x \sum_{x \ge d > \Delta}\frac{1}{d^2}
  = O\(\frac{x}{\Delta}\).
\end{split}
\end{equation}
Substituting~\eqref{eq:Sigma1 asymp} and~\eqref{eq:Sigma2 bound} in~\eqref{eq:split}, we obtain
\begin{equation}
\label{eq:Semi-final}
 T(x,y) = \Psi(x,y)\(\frac{6}{\pi^2}+O\(\frac{\log \log x \log \log y}{\log y} \)\)+O\(\frac{x}{\Delta}\).
\end{equation}

We now see  Lemmas~\ref{psi_rho} and~\ref{rho},  that for a sufficiently
 large $C$,  under the condition
$$
x \ge y \ge  \exp\(C \sqrt{\log x \log \log \log x}\),
$$
the bound~\eqref{psi_CEP} holds and furthermore,
we have
\begin{equation*}
\begin{split} \Psi(x,y)& ~ \frac{\log \log x \log \log y}{\log y}
\ge  \Psi(x,y) ~ \frac{1}{\log y}  \\
&\ge~ x \exp \(- 2 \frac{\log x}{\log y} \log  \frac{\log x}{\log y}
-\log \log y\) \\
& \ge ~ x \exp \(- 2 \frac{\log x}{\log y} \log \log x \)\\
& \ge ~ \max\left\{x \exp\(-\gamma \frac{ \log y \log \log y}{\log \log \log y}\), x^{1/2} (\log x)^\delta \right\}\\
& = ~ \frac{x}{\Delta}.
\end{split}
\end{equation*}
Therefore the term  $O\(x/\Delta\)$ can be removed
from~\eqref{eq:Semi-final}, which concludes the proof.
\end{proof}

\begin{thm}
\label{thm:Vxy} There  exists an absolute constant $C > 0$ such that for a sufficiently large $x$ the bound
$$
V(x,y)=
 \frac{3x}{\pi^2} + O\(\frac{ x\log \log x \log \log y}{\log y}\)
$$
holds uniformly  in the range
$$
x \ge y \ge  \exp\(C \sqrt{\log x \log \log \log x}\)
$$
where
$$
u=\frac{\log x}{\log y}.
$$
\end{thm}

\begin{proof}
Using partial summation, see~\cite[Chapter~I.0, Theorem~1]{Ten}, we can rewrite $V(x,y)$ in the following way,
\begin{eqnarray*}
 V(x,y)
&=&\frac{1}{\Psi(x,y)}  \sum_{\substack{a<n \le x \\
n\in S(x,y)\\ n\equiv a\pmod d}}\frac{\varphi(n - a)}{n - a}(n - a) \\
& =& \frac{1}{\Psi(x,y)}\(T(x,y) (x-a) - \int_1^x T(t,y) dt\).
\end{eqnarray*}
For
$$
  t \le x \qquad
 \text{and}
 \qquad  y \ge  \exp\(C \sqrt{\log x \log \log \log x}\)
$$
 Theorem~\ref{thm:Txy} implies that
$$
T(t,y)= \Psi(t,y) \(\frac{6}{\pi^2} + O\(\frac{\log \log t \log \log y}{\log y}\)\).
$$
 Therefore we get for $V(x,y)$
$$\frac{1}{\Psi(x,y)}\(x T(x,y)   - \(\frac{6}{\pi^2} + O\(\frac{\log \log x \log \log y}{\log y}\)\)
 \int_{1}^x \Psi(t,y) dt  \).
$$
Since $\Psi(t,y) \le \Psi(x,y)$, this simplifies as
\begin{equation}
\begin{split}
\label{eq:Prelim1}
 V(x,y)
&= ~ \frac{1}{\Psi(x,y)}\(x T(x,y)   - \frac{6}{\pi^2}
 \int_{1}^x \Psi(t,y) dt  \)\\
 &  \qquad \qquad \qquad \qquad
+~ O\( \frac{x\log \log x \log \log y}{\log y} \).
\end{split}
  \end{equation}
By Lemma~\ref{psi_rho} we have
\begin{eqnarray*}
\Psi(t,y) & = & t \rho\(\frac{\log t}{\log y}\)
\(1 + O\(\frac{ \log (u+1)}{\log y} \)\)\\
& = &t \rho\(\frac{\log t}{\log y}\)
+ O\(x \rho(u) \frac{ \log (u+1)}{\log y} \)\\
& = & t \rho\(\frac{\log t}{\log y}\) + O\(\Psi(x,y) \frac{ \log (u+1)}{\log y} \).
\end{eqnarray*}
Thus we derive from~\eqref{eq:Prelim1}
\begin{equation}
\label{eq:Prelim2} V(x,y) = \frac{1}{\Psi(x,y)}\(x T(x,y)   - \frac{6}{\pi^2} I(x,y) \)+O\(\frac{ x \log \log x \log \log y}{\log y} \),
\end{equation}
where
$$
 I(x,y) =
\int_{1}^x  t \rho\(\frac{\log t}{\log y}\) dt .
$$
Using integration by parts, we derive
\begin{eqnarray*}
 I(x,y) & = &
\frac{1}{2} \int_{1}^x    \rho\(\frac{\log t}{\log y}\) dt^2 \\
 & = & \frac{1}{2} x^2  \rho\(\frac{\log x}{\log y}\) + O(1)
- \frac{1}{2} \int_{1}^x  t^2 d \rho\(\frac{\log t}{\log y}\) \\
& = & \frac{1}{2} x^2  \rho\(u\) + O(1) - \frac{1}{2 \log y} \int_{1}^x  t  \rho'\(\frac{\log t}{\log y}\)  d t.
\end{eqnarray*}
By Lemma~\ref{rho'} we have
$$
\int_{1}^x  t  \rho'\(\frac{\log t}{\log y}\)  d t \ll \int_{1}^x  t  \rho\(\frac{\log t}{\log y}\)  d t \log (u+1) \ll I(x,y) \log (u+1).
$$
Therefore
$$
I(x,y) =  \frac{1}{2} x^2  \rho\(u\) + O\(1  + I(x,y) \frac{ \log (u+1)}{\log y}\)
$$
which, together with  Lemma~\ref{psi_rho}, implies
\begin{eqnarray*}
 I(x,y) & = & \frac{1}{2} x^2  \rho\(u\) \(1 +
O\( \frac{ \log (u+1)}{\log y}\)\) \\
 & = &
 \frac{1}{2} x \Psi(x,y) \(1 +
O\( \frac{ \log (u+1)}{\log y}\)\).
\end{eqnarray*}
Inserting this asymptotic formula in~\eqref{eq:Prelim2} and using Theorem~\ref{thm:Txy}, we conclude the proof.
\end{proof}

\section{Remarks}

Certainly improving the error term, or obtaining
similar bounds in a wider range are natural directions
for further investigation.

Studying average values of other number theoretic
functions on shifted smooth
numbers, such as
$$\frac{1}{\Psi(x,y)} \sum_{\substack{a<n \le x\\
n\in S(x,y)} }\tau(n - a) \qquad\text{and}\qquad
\frac{1}{\Psi(x,y)}\sum_{\substack{a<n \le x\\
n\in S(x,y)} }\omega(n - a),
$$
where $\tau(m)$ and $\omega(m)$ are the number
of positive integer divisors and the number
of prime divisors of $m \ge 1$, respectively, is
of ultimate interest too. However investigating
these sums may require a very different
approach.

\vspace{0.5cm}

\end{document}